\documentclass[a4paper,11pt]{article}
\usepackage[utf8x]{inputenc}
\usepackage{a4wide}
\usepackage{amssymb}
\usepackage{amsmath}
\usepackage{amsthm}
\usepackage{latexsym}
\usepackage[mathscr]{euscript}
\usepackage
{hyperref}
\usepackage[color]{changebar}
\usepackage{todonotes}
\usepackage{enumitem}

\usepackage{authblk}

\theoremstyle{definition}
\newtheorem{definition}{Definition}[section]

\theoremstyle{plain}
\newtheorem{Theorem}[definition]{Theorem}
\newtheorem{Proposition}[definition]{Proposition}
\newtheorem{Lemma}[definition]{Lemma}
\newtheorem{Corollary}[definition]{Corollary}

\theoremstyle{remark}
\newtheorem{remark}[definition]{Remark}
\newtheorem{exm}[definition]{Example}

\newcommand{\R}{\mathbb R}  
\newcommand{\N}{\mathbb N}

\newcommand{\eps}{\varepsilon}
\newcommand{\vphi}{\varphi}
\cbcolor{blue}

\newcommand{\sse}{\subseteq}

\newcommand{\tsf}{\rho} 

\newcommand{\LLSn}{Lorentzian length space}
\newcommand{\LpLS}{Lorentzian pre-length space }
\newcommand{\LpLSn}{Lorentzian pre-length space}
\newcommand{\Xll}{$(X,d,\ll,\leq,\rho)$ }

\newcommand{\Om}{\Omega}

\begin{document}

\title{Null distance and convergence of Lorentzian length spaces}

\author{Michael Kunzinger\thanks{E-mail: michael.kunzinger@univie.ac.at}}
\author{Roland Steinbauer\thanks{E-mail: roland.steinbauer@univie.ac.at}}

\affil{Faculty of Mathematics, University of Vienna,\\ Oskar-Morgenstern-Platz 1, A-1090 Wien, Austria.}

\date{\today}


\maketitle

\begin{abstract}
The null distance of Sormani and Vega encodes the manifold topology as well as the causality structure of a (smooth) spacetime. We extend this concept to Lorentzian length spaces, the analog of (metric) length spaces, which generalize Lorentzian causality theory beyond the manifold level. We then study Gromov-Hausdorff convergence based on the null  distance in warped product Lorentzian length spaces and prove first results on its compatibility with synthetic curvature bounds.
\vskip 1em

\medskip

\end{abstract}

\section{Introduction}

Metric geometry  \cite{BH99, BBI} 
has led to identifying the `metric core' of many results in Riemannian differential geometry, to 
clarifying the interdependence of various concepts, and to generalizations of central notions to lower regularity. In particular, Riemannian manifolds carry a natural metric structure and standard notions of convergence such as Gromov-Hausdorff (GH) and Sormani-Wenger intrinsic flat (SWIF) convergence \cite{SW11} interact well with geometric quantities, in particular with curvature bounds.

Despite the increasing demand for a Lorentzian analog of this framework, particularly driven by General Relativity (GR), see e.g.\ \cite{CGKE18}, a comparable metric theory is still in its infancy. Only recently Sormani and Vega \cite{SV16} put forward a solution to one prime obstacle i.e., the fact that the Lorentzian distance does not induce a metric structure\footnote{For a general discussion see \cite{V:21}.}. They constructed a `null distance' capable of encoding both the topological and the causality structure of the manifold, and first convergence results built on the corresponding metric and integral current structures were established by Burtscher and Allen \cite{AB19}.

In a somewhat parallel approach more directly rooted in GR, Kunzinger and S\"amann \cite{KS18} introduced a notion of Lorentzian length spaces. Based on the time separation function, this construction provides a close Lorentzian analog of (metric) length spaces. In particular, it allows one to extend beyond the manifold level the synthetic approach to (sectional) curvature bounds, which was introduced for general semi-Riemannian manifolds by Alexander and Bishop \cite{AB08}. 

In this work we provide a natural next step towards a comprehensive notion of metric limits for Lorentzian manifolds by considering GH convergence based on the null distance in the class of Lorentzian length spaces, thereby extending the results in \cite{AB19} beyond the manifold level. Furthermore, we show that in certain warped product Lorentzian length spaces GH-convergence interacts well with synthetic curvature bounds.  
\medskip

This work is structured in the following way: We collect preliminaries on the null distance, the convergence results of \cite{AB19}, Lorentzian synthetic curvature bounds and Lorentzian length spaces in Section \ref{sec:preliminaries}. Then, in Section \ref{sec:nulldist}, we extend the null distance to the setting of Lorentzian (pre-)length spaces and, following \cite{SV16,AB19}, establish its fundamental properties. Finally, in Section \ref{sec:warped}, we study GH limits of warped product Lorentzian length spaces and prove or main results on their interaction with curvature bounds.
\medskip

In the remainder of this introduction we collect some basic notions and conventions. All manifolds are assumed to be smooth, connected, Hausdorff, second countable, of arbitrary dimension $n\geq2$, and without boundary. A spacetime $(M,g)$ is a time oriented Lorentzian manifold, where we use the signature $(-,+,\dots,+)$. We will deal with metrics of various regularity, but generally assume them to be smooth unless explicitly stated otherwise. Causality notions will be based on locally Lipschitz curves and we denote the \emph{chronological} and the \emph{causal} relation by $I$ and $J$ and write $p\ll q$ and $p\leq q$ if $q\in I^+(p)$ and $p\in J^+(q)$, respectively. A \emph{generalized time function} $\tau:M\to \R$ is a function that is strictly increasing along all future directed causal curves. It is called a \emph{time function} if it is continuous. For points $p,q\in M$ the \emph{time separation function} $\rho(p,q)$ is the supremum of the length of all future directed causal curve segments from $p$ to $q$ with the understanding that $\rho(p,q)=0$ if there is no such curve, i.e., if $q\not\in J^+(p)$. In all matters of semi-Riemannian geometry and causality theory we will adopt the conventions and notations of \cite{ON83} and \cite{Min:19b}. 

Finally, we will often deal with \emph{warped products} $M=B\times_f F$, where $f>0$ is the warping function and $(B,g_B)$ and $(F,g_F)$ are semi-Riemannian manifolds. The metric $g$ on $M=B\times F$ is then given by  $g=g_B+f^2g_F$, where, as usual, we notationally suppress the projections. We will most of the time deal with the case that the base $B$ is a real interval $I$ and the fiber $(F,g_F)$ is Riemannian. Then the warped product metric takes the form $g=-dt^2+f^2g_F$.

\section{Preliminaries}\label{sec:preliminaries}

The null distance of Sormani and Vega \cite{SV16} provides a way of encoding the manifold topology as well as the causal structure of a spacetime. Given a generalized time function $\tau$ on $(M,g)$ the null distance is defined by
\begin{equation}
 \hat d_\tau(p,q)=\inf\{\hat L_\tau(\beta):\ \beta\ \mbox{piecewise causal from $p$ to $p$\}}\,.
\end{equation}
Here the null length of $\beta:[a,b]\to M$ is given by $\hat L_\tau(\beta)=\sum_{i=1}^k|\tau(x_i)-\tau(x_{i-1})|$, where
$a=s_0<s_1<\dots<s_k=b$, $x_i=\beta(s_i)$ are the break points.
It then holds that $\hat d_\tau$ is a pseudometric on $M$ \cite[Prop.\ 3.8]{SV16} and it is a (conformally invariant) metric that induces the manifold topology if $\tau$ is continuous and locally anti-Lipschitz, i.e.,
if any point in $M$ has a neighbourhood $U$ with a distance function $d_U$ such that  
\begin{equation}
 \forall\ x\leq y\ \in U\ \Rightarrow\ \tau(y)-\tau(x)\geq d_U(x,y).
\end{equation}
We say that $\hat d_\tau$ \emph{encodes the causality} of $M$ if 
\begin{equation}
 x\leq y \Leftrightarrow \hat d_\tau(x,y)=\tau (y)-\tau(x),
\end{equation}
a property that is stronger than definiteness \cite[Lem.\ 3.12]{SV16}.
For warped product spacetimes it holds by \cite[Thm.\ 3.25]{SV16} that the null distance induced by any smooth time function is definite and encodes the manifold topology as well as the causality. Note that the completeness assumption on the fiber is not needed, a fact which also follows from our Theorem \ref{th:d_t_definite_and_J}, below. 
\medskip

The metric and integral current structure of Lorentzian manifolds based on the null distance was further studied in \cite{AB19}. There, Allen and Burtscher showed that for any spacetime $(M,g)$ with locally anti-Lipschitz time function $\tau$, $(M,\hat d_\tau)$ is a length space \cite[Thm.\ 3.5]{AB19}. They also proved first GH and SWIF convergence results for warped product spacetimes: Given a (connected, compact) Riemannian manifold $(\Sigma,h)$, they consider sequences of warped product spacetimes $(M=I\times\Sigma,g_j=-dt^2+f_j^2(t)h)$ where $I$ is a closed interval. If the (continuous) warping functions $f_j\colon I\to (0,\infty)$ are uniformly bounded away from $0$ and if they converge uniformly to a limit function $f$, then the corresponding null distances $\hat d_j$ converge uniformly to $\hat d_f$ on $M=I\times\Sigma$ and the metric spaces $(M,\hat d_{g_j})$ converge to $(M,\hat d_{g}$) in the GH as well as in the SWIF sense \cite[Thm. 5.5]{AB19}. Here $g$ is the limiting warped product metric $g=-dt^2+f(t)^2h$.
\medskip

The Lorentzian length spaces of \cite{KS18} generalize the notion of length space to the Lorentzian world. Let $(X,d)$ be a metric space and assume $X$ is endowed with a preorder $\leq$ as well as a transitive relation $\ll$ contained in $\leq$, which we call the timelike and causal relation, respectively. If, in addition, we have a lower semicontinuous map\footnote{In previous accounts on Lorentzian length spaces the time separation function was denoted by $\tau$. To comply with our main points of reference \cite{SV16,AB19} we here reserve the letter $\tau$ for time functions.} $\rho \colon X\times X \to [0, \infty]$ that satisfies the reverse triangle inequality 
and $\rho(x,y)>0 \Leftrightarrow x\ll y$, then \Xll is called a \emph{Lorentzian pre-length space\/} with  \emph{time separation function\/} $\rho$.

A locally Lipschitz curve on an arbitrary interval $\gamma \colon I\rightarrow X$ that is non-constant on any sub-interval is called (future-directed) \emph{causal (timelike)} if for all $t_1<t_2\in I$ we have $\gamma(t_1)\leq\gamma(t_2)$ ($\gamma(t_1)\ll\gamma(t_2)$). A causal $\gamma$ is called \emph{null\/} if no two points on the curve are timelike related. 
The length of a future-directed causal $\gamma \colon [a,b]\rightarrow X$ is defined via $\rho$ by
\[
L_\rho(\gamma):=
\inf\Big\{\sum_{i=0}^{N-1} \rho(\gamma(t_i),\gamma(t_{i+1})): a=t_0<t_1<\ldots<t_N=b,\ N\in\N\Big\}.
\ 
\]
We call a future-directed causal curve $\gamma \colon [a,b]\rightarrow X$ 
\emph{maximal\/} if it realizes the time separation, $L_\rho(\gamma) = \rho(\gamma(a),\gamma(b))$. 
In analogy with metric length spaces we call $X$  a \emph{\LLSn\/} if, in addition to some 
technical assumptions (cf.\ \cite[Def. 3.22]{KS18}) $\rho = \mathcal{T}$, where for any  $x,y\in X$ we set
\begin{equation*}
\mathcal{T}(x,y):= \sup\{L_\rho(\gamma):\gamma \text{ future-directed causal from }x \text{ to } y\}\,, 
\end{equation*}
if there is a future-directed causal curve from $x$ to $y$. Otherwise we set $\mathcal{T}(x,y):=0$.

Causality theory in Lorentzian length spaces \cite{KS18,GKS19,HPS20} extends standard causality theory \cite{Min:19b} beyond the spacetime setting, to which it reduces for smooth strongly causal spacetimes. Hence any smooth strongly causal spacetime is an example of a Lorentz\-ian length space, but more generally,
spacetimes with low regularity metrics and certain Lorentz-Finsler spaces \cite{Min19} provide further examples \cite[Sec.\ 5]{KS18}. In particular, any continuous spacetime with strongly causal and causally plain metric (a condition that rules out causal pathologies, see \cite[Def.\ 1.16]{CG}) is a (strongly localizable, for a definition see below) Lorentzian length space.
\medskip

Based on pioneering work by Harris \cite{H82}, Alexander and Bishop in \cite{AB08} gave a characterization of sectional curvature bounds in terms of triangle comparison in smooth semi-Riemannian manifolds. We say that $(M,g)$ has \emph{sectional curvature bounded below} by some constant $K$, $R\geq K$, if the sectional curvatures for all spacelike planes are bounded below by $K$ and if for all timelike planes they are bounded above by $K$. Equivalently we have 
\begin{equation}
 R\ge K\quad\mbox{if}\quad  
 R(v,w,v,w)\geq K\left(\langle v,v,\rangle\langle w,w,\rangle-\langle v,w\rangle^2\right).
\end{equation}
Then it holds \cite[Thm.\ 1.1]{AB08} that $R\geq K$ ($R\leq K$) if and only if in any convex (totally normal) neighbourhood the signed length of the geodesic between two points on a geodesic triangle is at least (at most) that of the corresponding points in the model triangle in $\mathbb{L}^2(K)$. Here the \emph{signed length} of a geodesic in a convex neighbourhood is defined as the signed length of the connecting vector $|\gamma_{pq}|_\pm=\text{sign}(\gamma_{pq})\sqrt{|\langle \gamma_{pq},\gamma_{pq}\rangle|}$, with the sign of timelike vectors taken to be negative. Moreover, the Lorentzian model spaces $\mathbb{L}^2(K)$ of constant sectional curvature $K$ are 
\begin{equation}\label{eq:model_spaces}
\mathbb{L}^2(K) = \left\{ \begin{array}{ll}
\tilde S^2_1(r) & K=\frac{1}{r^2}\\
\R^2_1 & K=0\\
\tilde H^2_1(r) & K= -\frac{1}{r^2}\,,
\end{array}
\right.
\end{equation}
where  $\tilde S^2_1(r)$ is the simply connected covering manifold of the two-dimensional \emph{Lorentzian pseudosphere} $S^2_1(r)$, $\R^2_1$ is two-dimensional \emph{Minkowski space}, and $\tilde H^2_1(r)$ is the simply connected covering manifold of the two-dimensional \emph{Lorentzian pseudohyperbolic} space.
\medskip

Again, in parallel to the case of metric geometry, appropriate notions of synthetic (timelike or causal) curvature bounds based on triangle comparison have been introduced in Lorentzian length spaces.
By a {timelike geodesic triangle\/} we mean a triple $(x,y,z)\in X^3$ with $x\ll y \ll z$ such that $\rho(x,z)<\infty$ and such that the sides are realized by future-directed causal curves. 
We then say, cf.\ \cite[Def.\ 4.7]{KS18}, that a Lorentzian pre-length space \Xll has \emph{timelike curvature bounded below} (above) by 
$K\in\R$ if every point in $X$ has a so-called comparison neighborhood $U$ such that:
\begin{enumerate}[label=(\roman*), topsep=2pt,itemsep=2pt,parsep=0pt,partopsep=0pt]
	\item $\rho|_{U\times U}$ is finite and continuous.
	\item Whenever $x,y \in U$ with $x \ll y$, there exists a causal curve $\alpha$ in $U$ with $L_\rho(\alpha) = 
\rho(x,y)$.
	\item If $(x,y,z)$ is a timelike geodesic triangle in $U$, realized by maximal causal curves $\alpha, 
\beta, \gamma$ whose side lengths satisfy the appropriate size restrictions (see \cite[Lem.\ 4.6]{KS18}), 
and if $(x',y',z')$ is a comparison triangle of $(x,y,z)$ in $\mathbb{L}^2(K)$ realized by timelike geodesics $\alpha '$, $\beta '$, $\gamma '$, then whenever $p$, $q$ are points on the sides of $(x,y,z)$ and $p', q'$ 
are corresponding points\footnote{This means that $p'$ lies on the side corresponding to the side containing $p$ at the same time separation from the vertex (i.e., e.g.\ if $p$ lies on the side $xy$ then $\rho(x,p)=\rho'(x',p')$, etc.). Similarly for $q'$.} of $(x',y',z')$, we have $\rho(p,q)\le \rho '(p',q')$ 
$($respectively $\rho(p,q)\ge \rho '(p',q'))$.
\end{enumerate}
\medskip

We close this section by recalling further central notions in Lorentzian pre-length spaces. Generally, causality conditions such as strong causality and global hyperbolicity are translated in perfect analogy from the spacetime setting. \Xll is called \emph{causally path connected} if for all $x,y\in X$ with $x\ll y$ and all $x,y$ with $x\leq y$ there is a future-directed timelike resp.\ causal curve from $x$ to $y$. 
A \emph{localizing neighborhood} $\Omega_x$ of a point $x\in X$ is a substitute for a convex neighbourhood, and a Lorentzian pre-length space is called \emph{localizable} if every $x$ has such a neighbourhood. It is defined by the conditions:
\begin{enumerate}[label=(\roman*), topsep=2pt,itemsep=2pt,parsep=0pt,partopsep=0pt]
\item \label{def-loc-LpLS-cau-comp} There is a $C>0$ such that $L^d(\gamma)\leq C$ for all causal curves $\gamma$ contained in $\Omega_x$ (we say that $X$ is $d$-compatible).
  \item \label{def-loc-LpLS-om-con} There is a continuous map $\omega_x\colon \Omega_x \times \Omega_x\rightarrow [0,\infty)$ such that
$(\Omega_x, d\rvert_{\Omega_x\times\Omega_x}, \ll\rvert_{\Omega_x\times \Omega_x}, \leq\rvert_{\Omega_x\times\Omega_x}, 
\omega_x)$ is a \LpLS with the following non-triviality condition: for every $y\in\Omega_x$ we have 
$I^\pm(y)\cap\Omega_x\neq\emptyset$.
  \item \label{def-loc-LpLS-max-cc} For all $p,q\in \Omega_x$ with $p<q$ there is a future-directed causal curve 
$\gamma_{p,q}$ from $p$ to $q$ that is 
maximal in $\Omega_x$ and satisfies
\begin{equation}
 L_\tau(\gamma_{p,q}) = \omega_x(p,q) \leq \tau(p,q)\,.
\end{equation}
\end{enumerate}

If, in addition, the neighborhoods $\Omega_x$ can be chosen such that
\begin{enumerate}[label=(\roman*), topsep=2pt,itemsep=2pt,parsep=0pt,partopsep=0pt]
 \item[(iv)]\label{def-loc-LpLS-4} Whenever $p,q\in\Omega_x$ satisfy $p\ll q$ then $\gamma_{p,q}$ is timelike and strictly longer than any future-directed 
causal curve in $\Omega_x$ from $p$ to $q$ that contains a null segment,
\end{enumerate}
then \Xll is called {\em regularly localizable}. Finally, if every point $x\in X$ has a neighborhood basis
of open sets $\Omega_x$ satisfying (i)--(iii), respectively (i)--(iv), then \Xll is called {\em strongly localizable} respectively {\em SR-localizable}. 

In a strongly causal and localizable Lorentzian pre-length space the length $L_\rho$ is upper semicontinuous, if it is regularly localizable, maximal causal curves have a causal character and the push-up principle holds \cite[Prop.\ 3.17, Thms.\ 3.18, 3.20]{KS18}. A \LpLS is called \emph{geodesic} if for all $x<y$ there is a future-directed causal curve $\gamma$ from $x$ to $y$ with $\tau(x,y)=L_\tau(\gamma)$ (hence maximizing). Any globally hyperbolic \LpLS is geodesic \cite[Thm.\ 3.30]{KS18}.

\section{The null distance in Lorentzian length spaces}\label{sec:nulldist}
In this section we extend the notion of null distance to the setting of Lorentzian (pre-)length spaces and establish its fundamental properties. 

\begin{definition}
Let $(X,d,\ll,\le,\tsf)$ be a \LpLSn. A map $\tau: X\to \R$ is called a \emph{generalized time function} if $\tau$ is strictly
increasing along every (non-trivial) future-directed causal curve. If $\tau$ is continuous, it is called a \emph{time function}.
\end{definition}

As we shall see below in Thm.\ \ref{th:top_anti_lip} existence of a ``reasonable'' time function implies strong causality. While a significant part of the causal ladder for Lorentzian length spaces has been established in \cite{HPS20,KS18}, the precise relationship between existence of time functions and stable causality in this general framework is still an open question.

\begin{definition}
A map $\beta: [a,b] \to X$ from a closed interval into a \LpLS $X$ is called a \emph{piecewise causal curve} if there exists a 
partition $a=s_0<s_1<\dots<s_{k-1}<s_k=b$ such that each $\beta_i:=\beta|_{[s_{i-1},s_i]}$ is either trivial (i.e., constant) or 
future directed causal or past directed causal. 
Given, in addition, a generalized time function $\tau: X\to \R$, the \emph{null length} of $\beta$ is
\begin{equation}\label{def:null_length}
\hat L_\tau(\beta) := \sum_{i=1}^k |\tau(x_i)-\tau(x_{i-1})|,
\end{equation}
where $x_i=\beta(s_i)$ ($i=0,\dots,k$). Moreover, for $p, q\in X$ we define the \emph{null distance} of $p$ and $q$ by
\begin{equation}\label{def:null_distance}
\hat d_\tau(p,q) := \inf \{\hat L_\tau(\beta) \mid \beta \text{ piecewise causal from } p \text{ to } q \},
\end{equation} 
\end{definition}

\begin{remark}\label{rem:causal_non_constant} Contrary to the convention used in \cite[Def.\ 3.1]{SV16}, causal curves in \LpLSn s are always assumed to be
nowhere constant (in accordance with the common custom in general relativity). To obtain a faithful generalization of null length
and null distance from \cite{SV16} to our setting we therefore explicitly allowed constant (sub-)curves in the definition of
piecewise causal curves above.
\end{remark}
\begin{definition}
A \LpLS is called \emph{sufficiently causally connected} (scc) if it is path connected, causally path connected and if every point $p\in X$
lies on some timelike curve.
\end{definition}
Under this assumption we indeed have the following fundamental existence result:
\begin{Lemma}\label{lem:existence_of_pw_causal_curves} Let $(X,d,\ll,\le,\tsf)$ be an scc \LpLSn. 
Then for any $p, q\in X$ there exists a piecewise causal curve from $p$ to $q$.
\end{Lemma}
\begin{proof}  
By \cite[Lemma 2.12]{KS18}, each $I^\pm(x)$ ($x\in X$) is open, and by causal path connectedness the relation
$p\ll q$ is always realized by the existence of a future directed timelike curve from $p$ to $q$. 
Moreover, since any $p\in X$ lies on a timelike curve we have $X=\bigcup_{x\in X} I^-(x) \cup \bigcup_{y\in X} I^+(y)$. 
Based on these observations, a straightforward adaptation of the proof of \cite[Lemma 3.5]{SV16} yields the claim.
The only difference is that in the present situation the finite covering of any path from $p$ to $q$ will 
in general contain both timelike futures and timelike pasts of points in $X$\footnote{Note that this can really occur, e.g.\ for $[a,b]\times_f X$ in points with $t=a,b$.}.
\end{proof}
We also have the following analogue of \cite[Lemma 3.6]{SV16}:
\begin{Lemma}\label{lem:null_length_basics} Let $\tau$ be a generalized time function on a \LpLS and let $\beta: [a,b]\to X$ be piecewise
causal from $p$ to $q$. Then
\begin{itemize}
\item[(i)] $\hat L_\tau(\beta) = 0$ if and only if $\beta$ is trivial.
\item[(ii)] $\hat L_\tau(\beta) = \tau(q) - \tau(p)$ if and only if $\beta$ is future directed causal or constant.
\item[(iii)] $\hat L_\tau(\beta) \ge \max_{y\in\beta} \tau(y) - \min_{x\in \beta}\tau(x) \ge |\tau(q)-\tau(p)|$.
\item[(iv)] If $\tau\circ \beta: [a,b]\to X$ is absolutely continuous, then
\[
\hat L_\tau(\beta) = \int_a^b |(\tau\circ \beta)'|(s)\,ds.
\]
\end{itemize}
\end{Lemma}
\begin{proof} By \cite[Def.\ 2.18]{KS18} (cf.\ Remark \ref{rem:causal_non_constant}), causal
curves are always assumed to be non-constant, so (i) follows. The other properties are direct consequences of the definitions. 
\end{proof}
The following result collects basic properties of the null distance (cf.\ \cite[Lemma 3.8]{SV16}).
\begin{Lemma} Let $\tau$ be a generalized time function on an scc \LpLSn. Then the null distance $\hat d_\tau$ is a 
finite pseudometric.
\end{Lemma}
\begin{proof} Symmetry and triangle inequality are immediate from the definition, and finiteness follows from Lemma \ref{lem:existence_of_pw_causal_curves}.
That $\hat d_\tau(p,p)= 0$ is seen by considering the constant (hence piecewise causal) curve $\beta \equiv p$.
\end{proof}
Since the previous proof relied on Lemma \ref{lem:existence_of_pw_causal_curves}, we will henceforth usually assume the scc property
in order to avoid degenerate situations. The next result corresponds to Lemmas 3.10--3.13 and 3.16--3.18 in \cite{SV16} (with identical proofs):
\begin{Proposition}\label{prop:gen_tf} Let $\tau$ be a generalized time function on an scc \LpLS $X$. Then
\begin{itemize}
\item[(i)] For any $p, q\in X$, $\hat d_\tau(p,q)\ge |\tau(q)-\tau(p)|$. In particular, $\hat d_\tau(p,q)=0$ implies $\tau(p)=\tau(q)$.
\item[(ii)] If  $p\le q$, then $\hat d_\tau(p,q) = \tau(q)-\tau(p)$.
\item[(iii)] $\tau$ is bounded on causal diamonds: $p\le x\le q \Rightarrow \tau(p)\le \tau(x)\le \tau(q)$.
\item[(iv)] $\hat d_\tau$ is bounded on causal diamonds: $p\le x, y\le q \Rightarrow \hat d_\tau(x,y)\le 2(\tau(q) -\tau(p))$.
\item[(v)] If $X$ has the property that $p\le q \Leftrightarrow \hat d_\tau(p,q) = \tau(q)-\tau(p)$, then $\hat d_\tau$ is definite.
\item[(vi)] If $\tilde\tau$ is another generalized time function on $X$ and $\lambda\in (0,\infty)$, 
then $\hat d_{\tilde\tau}= \lambda\hat d_{\tau} \Leftrightarrow \tilde\tau = \lambda \tau + C$ for some constant $C$.
\end{itemize}
\end{Proposition}
Inspection of the proof of \cite[Prop.\ 3.14]{SV16} shows that, since timelike futures and pasts are open in any Lorentzian
pre-length space (\cite[Lemma 2.12]{KS18}), it carries over unchanged to the present setting. Therefore we have:
\begin{Proposition}\label{3.14} Let $\tau$ be a generalized time function on an scc \LpLS $X$. Then the following are equivalent:
\begin{itemize}
\item[(i)] $\tau: X\to \R$ is continuous.
\item[(ii)] $\hat d_\tau: X\times X \to \R$ is continuous.
\end{itemize}
\end{Proposition}
The following result is a direct generalization of \cite[Prop.\ 3.15]{SV16}. For the reader's convenience we adapt its
proof to the current, topologically more general, setup.
\begin{Proposition}\label{prop:hat_d_induces_top}
Let $(X,d,\ll,\le,\tsf)$ be an scc Lorentzian pre-length space with generalized time function $\tau$ and suppose that $(X,d)$ 
is locally compact. Then the following are equivalent:
\begin{itemize}
\item[(i)] $\hat d_\tau$ induces the same topology as $d$.
\item[(ii)] $\tau$ is continuous and $\hat d_\tau$ is definite.
\end{itemize}
\end{Proposition}
\begin{proof}
(i)$\Rightarrow$(ii): By Proposition \ref{3.14}, $\tau$ is continuous. Also, since $d$ is definite, the topology on $X$
is Hausdorff, implying that $\hat d_\tau$ is definite as well. 

(ii)$\Rightarrow$(i): Also in this case, $\hat d_\tau$ is continuous by Proposition \ref{3.14}, and so it only remains toshow that the $\hat d_\tau$-topology $\mathcal{O}_\tau$ is finer than the $d$-topology $\mathcal{O}_d$. 
Let $x\in X$ and $\eps>0$ and denote by  $B_\eps^d(x)$ the open $\eps$-ball for $d$. 
Since $X$ is locally compact we may assume $\eps$ small enough so that $\partial B_\eps^d(x)$ is compact.
We now distinguish two cases:
First, if $\partial B_\eps^d(x) = \emptyset$, then $X$ being connected implies $B_\eps^d(x)=X$, and so we 
may pick any $\eps_0>0$ to obtain $B_{\eps_0}^{\hat d_\tau}(x) \subseteq B_\eps^d(x)$.
Otherwise, $\eps_0:= \min \{\hat d_\tau(x,y) \mid y\in \partial B_\eps^d(x) \}$ exists and is 
strictly positive since $\hat d_\tau$ is definite. Let $y\not\in B_\eps^d(x)$ and pick a piecewise causal curve
$\beta:[a,b]\to M$ from $x$ to $y$. With $z$ the first intersection of $\beta$ with $\partial B_\eps^d(x)$, let 
$\beta_0$ be the initial part of $\beta$ from $x$ to $z$. Then $\hat L_\tau(\beta) \ge \hat L_\tau(\beta_0)\ge \eps_0$.
Thus also $\hat d_\tau(x,y)\ge \eps_0$. We conclude that also in this case $B_{\eps_0}^{\hat d_\tau}(x) \subseteq B_\eps^d(x)$,
so indeed $\mathcal{O}_\tau \supseteq \mathcal{O}_d$.
\end{proof}
Our next goal is to establish that the null distance is definite if and only if the generalized time function is 
anti-Lipschitz. First, for locally compact $X$ the following result (\cite[Lemma 4.3]{SV16}) carries over, with identical proof:
\begin{Lemma}\label{lem:def_d_tau}
Let $(X,d,\ll,\le,\tsf)$ be an scc Lorentzian pre-length space with generalized time function $\tau$ such that $(X,d)$ 
is locally compact. Let $d_U$ be a definite distance function on the open subset  $U\sse X$ such that
for all $x,y\in U$ we have $x\le y \Rightarrow \tau(y)-\tau(x)\ge d_U(x,y)$. Then for any $p\in U$ and any
$q\in X\setminus \{p\}$, $\hat d_\tau(p,q)>0$. In particular, $\hat d_\tau$ is definite on $U$.
\end{Lemma}
Following \cite[Def.\ 4.4]{SV16} we call a map $f: X\to \R$ \emph{anti-Lipschitz} on $U\sse X$ if there exists
a definite distance function $d_U$ on $U$ such that for all $x\le y$ in $U$ we have $f(y)-f(x)\ge d_U(x,y)$.
It is called \emph{locally anti-Lipschitz} if every point in $X$ possesses a neighborhood on which $f$ is
anti-Lipschitz. Such a function then is automatically a generalized time function, and by Lemma \ref{lem:def_d_tau} the 
corresponding null distance is definite. Together with Proposition \ref{prop:gen_tf} (ii) we obtain
(cf.\ \cite[Prop.\ 4.5]{SV16}):
\begin{Proposition}\label{prop:hat_d_def_tau_antiLip}
Let $\tau$ be a generalized time function on a locally compact ssc Lorentzian pre-length space $X$. Then $\hat d_\tau$ is definite
if and only if $\tau$ is locally anti-Lipschitz.
\end{Proposition}
In the above definition of the anti-Lipschitz property there is no requirement on the compatibility of the local 
distance function $d_U$ with the topology on $X$ induced by the metric $d$. To introduce such an assumption, we call
$f: X\to \R$ \emph{topologically anti-Lipschitz} on the open set $U\sse X$ if there exists a  
definite distance function $d_U$ on $U$ that induces the $d$-topology on $U$ and
such that for all $x,y\in U$ we have $x\le y \Rightarrow \tau(y)-\tau(x)\ge d_U(x,y)$.
The function is called \emph{topologically locally anti-Lipschitz} if every point in $X$ is contained in an open set $U$
equipped with such a $d_U$. Of course the simplest (and most relevant) situation is the one where indeed $d_U=d$.
As the following result shows, it is the topological anti-Lipschitz property that allows us to 
position the existence of time functions on the causal ladder of Lorentzian pre-length spaces.

\begin{Theorem}\label{th:top_anti_lip} Let $(X,d,\ll,\le,\tsf)$ be an scc \LpLS that is equipped with a topologically locally anti-Lipschitz
time function $\tau$. Then for each point $p\in X$ and each neighborhood $V$ of $p$ there exists a neighborhood
$U\sse V$ of $p$ that is causally convex: if $\gamma: [a,b]\to U$ is causal and $\gamma(a), \gamma(b)\in U$ then
$\gamma([a,b])\sse U$. In particular, if $X$ is a Lorentzian length space, then $X$ is strongly causal.
\end{Theorem}

\begin{proof} Fixing $p\in X$ and a neighborhood $V$ of $p$, pick $U$ and $d_U$ as in the assumption.
Denote the closed (resp.\ open) $d$-ball $\{q \mid d(p,q)\le \delta\}$ (resp.\ $\{q \mid d(p,q) < \delta\}$) 
of radius $\delta$ by $D^d_\delta(p)$ (resp.\ $B^d_\delta(p)$), and 
analogously for $d_U$ and pick $\delta>0$ such that $D^d_{2\delta}(p)\sse U$. Next, let $\eps>0$ be such that
$D^{d_U}_{2\eps}(p)\sse B^d_\delta(p)$ and define $\vphi_\eps: X\to \R$ by 
\[
\vphi_\eps(q) := \left\{
\begin{array}{ll}
\eps - \frac{1}{2} d_U(p,q) & q\in D^{d_U}_{2\eps}(p)\\
0 & \text{otherwise.}
\end{array}
\right.
\]
Then $\vphi_\eps$ is continuous and we claim that $\tau_\eps:= \tau + \vphi_\eps$ and $\tau_{-\eps}:= \tau - \vphi_\eps$ are time 
functions on $X$. Let us verify this for $\tau_\eps$, the proof for $\tau_{-\eps}$ being analogous. If $q_1, q_2\in U$ with 
$q_1 > q_2$, then noting that $|\vphi_\eps(q_1)-\vphi_\eps(q_2)|\le \frac{1}{2} d_U(q_1,q_2)< d_U(q_1,q_2)$ due to the reverse triangle 
inequality for $d_U$, we get
\begin{align*}
\tau_\eps(q_1) = \tau(q_1) + \vphi_\eps(q_1) \ge \tau(q_2) + d_U(q_1,q_2) + \vphi_\eps(q_1) > \tau_\eps(q_2).
\end{align*}
Now suppose that $q_1> q_2$ with $q_1\not\in U$ and $q_2\in U$. Then $\tau_\eps(q_1)=\tau(q_1)$ and either $\tau_\eps(q_2)=\tau(q_2)$, 
in which case we are done since $\tau$ is a time function, or $q_2\in D^{d_U}_{2\eps}(p)$. In the latter case, let $\gamma$ be a future 
directed causal curve from $q_2$ to $q_1$ and let $\bar q$ be an intersection of $\gamma$ with the boundary of $D^{d_U}_{2\eps}(p)$. Then
\[
\tau_\eps(q_1) = \tau(q_1) \ge \tau(\bar q) = \tau_\eps(\bar q) > \tau_\eps(q_2).
\]
The case $q_1\in U$, $q_2\not\in U$ follows symmetrically.

Now set $U_\eps(p) := \{q\in X \mid \tau_{-\eps}(q) < \tau(p) < \tau_\eps(q) \}$. Then $U_\eps(p)$ is an open neighborhood of $p$
that is contained in $U$ and we are going to show that $U_\eps(p)$ is causally convex. To this end, first note that
\begin{equation}\label{eq:boundary}
\partial U_\eps(p) \sse \{q \mid \tau_\eps(q) = \tau(p) \} \cup \{q \mid \tau_{-\eps}(q) = \tau(p) \} =: \Gamma_\eps \cup \Gamma_{-\eps}.
\end{equation}
If $\gamma$ is a future directed causal curve that intersects $\Gamma_\eps$ at $\gamma(t_0)$ then since $\tau_\eps$
is a time function it follows that $\tau_\eps(\gamma(t))\le \tau_\eps(\gamma(t_0))=\tau(p)$ for $t\le t_0$, so in particular
$\gamma(t)\not\in U_\eps$ for $t\le t_0$.
Analogously, if $\gamma$ intersects $\Gamma_{-\eps}$ at $\gamma(t_0)$ then $\gamma(t)\not\in U_\eps$ for $t\ge t_0$.
Together with \eqref{eq:boundary} this implies that $\gamma$ can enter $U_\eps(p)$ only by
passing through $\Gamma_\eps$ and can leave it only through $\Gamma_{-\eps}$. Suppose now that $\gamma: [a,b] \to X$ is future directed causal
with $\gamma(a), \gamma(b) \in U_\eps(p)$. If $\gamma$ were not entirely contained in $U_\eps(p)$ then by what was just shown there would exist $s<t$ in $(a,b)$ such that
$\tau_{-\eps}(\gamma(s)) = \tau(p) =\tau_\eps(\gamma(t))$. But then $\tau_\eps(\gamma(s))\ge \tau_{-\eps}(\gamma(s)) = \tau_\eps(\gamma(t))$,
contradicting the fact that $\tau$ is a time function.

The final claim follows from \cite[Th.\ 3.26 (iv)]{KS18}.
\end{proof}

\begin{exm} Let $(M,g)$ be a smooth manifold with a continuous, causally plain Lorentzian metric $g$ and an arbitrary background Riemannian metric $h$. Alternatively, we may consider a locally Lipschitz proper Lorentz-Finsler space $(M,\mathcal{F})$
in the sense of \cite{Min19} such that for the corresponding cone structure $C$ we have $\mathcal{F}(\partial C)=0$.
In both cases we then obtain a Lorentzian pre-length space $(M,d^h,\ll,\le,\rho)$ (see \cite[Sec.\ 5]{KS18}). Assume that
$\tau$ is an $h$-steep time function on $M$ (see \cite[Sec.\ 2.2]{Min19}). Then if $\gamma:[a,b]\to M$ is a future directed causal curve
from $p$ to $q$ we have
\[
\tau(q)-\tau(p) = \int_a^b d\tau(\dot\gamma)(t)\,dt \ge \int_a^b \|\dot\gamma(t)\|_h\,dt \ge d_h(p,q).
\]
Thus $\tau$ is topologically locally anti-Lipschitz and Theorem \ref{th:top_anti_lip} implies that any point in $M$
has a neighborhood base consisting of causally convex sets. In both cases this implies that $(M,d^h,\ll,\le,\rho)$ is 
a Lorentzian length space (cf.\ \cite[Th.\ 5.12, 5.16]{KS18}). 
Of course for smooth Lorentzian manifolds and even for the much more general
Lorentz-Finsler setting this is well known. Indeed, for closed cone structures the existence of a time function implies 
stable causality, which in turn implies strong causality (\cite[Th.\ 2.30]{Min19}). It is not known whether an
analogous result also holds for Lorentzian length spaces (cf.\ \cite{HPS20} for the definition of stable causality of 
Lorentzian length spaces).
\end{exm}
\begin{definition}
A piecewise causal curve $\beta: [a,b]\to X$ in an scc Lorentzian pre-length space $X$ with generalized time function $\tau$
is called \emph{minimal} if it minimizes the null distance, i.e., if $\hat d_\tau(\beta(a),\beta(b))=\hat L_\tau(\beta)$.
\end{definition}
As in \cite[Cor.\ 3.19]{SV16} it follows directly from Lemma \ref{lem:null_length_basics} (ii) that any causal curve $\beta$
from $p$ to $q\ge p$ is minimal with
\begin{equation}\label{eq:causal_min}
\hat d_\tau(p,q) = \hat L_\tau(\beta) = \tau(q) - \tau(p).
\end{equation}
In other words, causal curves are null distance realizers. 

Next, we transfer \cite[Lemma 3.20]{SV16} to the Lorentzian length space setting:
\begin{Proposition}
Let $\beta: [a,b]\to X$ be a minimal piecewise causal curve in an scc Lorentzian pre-length space $X$ with generalized time 
function $\tau$. Then either $\beta$ is causal or it is piecewise null, changing direction at each break point. If $X$ is, 
in addition, localizable, then $\beta$ is a piecewise null geodesic.
\end{Proposition}
\begin{proof}
The proof of \cite[Lemma 3.20]{SV16}, relying only on the properties that $I^\pm(p)$ is open for any $p\in X$ and that
$\tau$ is strictly increasing along future directed causal curves, shows that $\beta$ has the following 
property: If $a<s_0<s_1<\dots < s_k =b$ denotes the breaks of $\beta$, then no two points of $\beta|_{[s_{i},s_{i+1}]}$
are timelike related, giving the first claim (cf.\ \cite[Def.\ 2.18]{KS18}). Suppose now that $X$ is localizable, let
$t_0\in (s_i,s_{i+1})$, $x_0:=\beta(t_0)$ and let $\Om_{x_0}$ be a localizing neighborhood of $x_0$. Pick $\eps>0$ such that
$\beta([t_0-\eps,t_0+\eps])\sse \Om_{x_0}$. If $\beta|_{[t_0-\eps,t_0+\eps]}$ were not maximizing for $\rho$, then
there would exist a causal curve $\gamma_{p,q}$ from $p=\beta(t_0-\eps)$ to $q=\beta(t_0+\eps)$ in $\Om_{x_0}$
such that $L_\rho(\beta|_{[t_0-\eps,t_0+\eps]}) < L_\rho(\gamma_{p,q})$. But then, in particular, $0< L_\rho(\gamma_{p,q})$
and thereby $p\ll q$, a contradiction. Consequently, $\beta|_{[s_{i},s_{i+1}]}$ is a geodesic in the sense of 
\cite[Def.\ 4.1]{GKS19}.
\end{proof}
Finally, we immediately conclude the analogue of \cite[Cor.\ 3.21]{SV16}:
\begin{Corollary} Let $X$ be an scc Lorentzian pre-length space $X$ with generalized time 
function $\tau$ and suppose that $p, q\in X$ are not causally related. If $\beta$ is a piecewise causal curve 
from $p$ to $q$ that contains a timelike subsegment then there exists a strictly shorter piecewise causal curve 
$\alpha$ from $p$ to $q$, i.e., with $\hat L_\tau(\alpha) < \hat L_\tau(\beta)$.
\end{Corollary}

\section{Warped products}\label{sec:warped}

Warped products are of fundamental importance in Riemannian and Lorentzian geometry. In the context of general relativity they are generalizations of Friedmann-Lema\^{\i}tre-Robertson-Walker spacetimes, which serve as basic  cosmological models of our universe. Warped products and generalized cones likewise play a fundamental role in length spaces with synthetic curvature bounds. In the context of Lorentzian length spaces they have been studied in \cite{AGKS20} and we refer to this work for further background information. We start with a closer look at the null distance on such spaces.

\subsection{Null distance on generalized cones}

To begin with, we recall some basics of warped products with one-dimensional basis, the so-called \emph{generalized cones} from \cite[Sec.\ 3,4]{AGKS20}.

Let $(X,d)$ be a metric space and $I\subseteq\R$ an interval. Observe that in \cite{AGKS20} only the case of open $I$ has been investigated. Here, to also allow for the compact case, we keep $I$ general and often write $I=\langle a,b\rangle$ if wee need to specify the boundary points. 
Then set $Y:=I\times X$ and put the product metric on  $Y$, i.e., $D((t,x),(t',x')) = |t-t'| + d(x,x')$ for $(t,x),(t',x')\in Y$. Let $f\colon I\rightarrow (0,\infty)$ be 
continuous.  Let $\gamma\colon J\rightarrow Y$, $\gamma=(\alpha,\beta)$, where $\alpha\colon J\rightarrow I$ and $\beta\colon 
J\rightarrow X$ are both absolutely continuous and $\alpha$ is strictly monotonous. Then the metric derivative $v_\beta$ of the curve $\beta$ exists almost everywhere \cite[Thm.\ 1.1.2]{AGS:05} and we call $\gamma$
\begin{equation}\label{eq:tl_null_causal}
\begin{cases}
 \text{timelike}\\
 \text{null}\\
 \text{causal}
 \end{cases}
 \text{\quad if \qquad}
 -\dot\alpha^2 + (f\circ \alpha)^2 v_{\beta}^2\quad
 \begin{cases}
  < 0\\
  = 0\\
  \leq 0\,,
 \end{cases}
\end{equation}
almost everywhere. It is called \emph{future/past directed causal} if $\alpha$ is strictly monotonically 
increasing/decreasing, i.e., $\dot\alpha>0$ or $\dot\alpha<0$ almost everywhere. The length $L(\gamma)$ of a causal curve $\gamma$ is defined by
\begin{equation*}
 L(\gamma):= \int_a^b \sqrt{\dot\alpha^2 - (f\circ\alpha)^2 v_\beta^2}\,.
\end{equation*}
The time separation function $\rho\colon Y\times Y \rightarrow [0,\infty]$ (called $\tau$ in \cite{AGKS20}) is defined as
\begin{equation}
 \rho(y,y'):=\sup\{L(\gamma):\gamma \text{ future directed causal curve from } y \text{ to } y'\}\,,
\end{equation}
if this set is non-empty, and $\rho(y,y'):= 0$ otherwise. The causal and timelike relations $p\le q$ resp.\ $p\ll q$ are defined as usual via the existence of a causal resp.\ timelike future directed curve from $p$ to $q$. When $I\times X$ is equipped with this time separation function,
we write $Y\equiv I\times_fX$ and call it a \emph{generalized cone} (or warped product with one-dimensional basis) with \emph{warping 
function} $f$.

According to \cite[Prop.\ 3.26]{AGKS20}, if $Y=I\times_f X$ is a generalized cone, where $(X,d)$ is a  
length space,  then $(Y,D,\ll,\leq,\rho)$ is a Lorentzian pre-length space. 
By \cite[Th.\ 4.8, Cor.\ 4.9]{AGKS20} we have:
\begin{Theorem}
Any generalized cone $I\times_f X$, where $I$ is open and $(X,d)$ is a locally compact length space, is a strongly causal Lorentzian length space. If $X$ is, in addition, geodesic, then $I\times_f X$ is a 
regular\footnote{See \cite[Def.\ 3.22]{KS18} for a definition.} Lorentzian length space.
\end{Theorem}
Furthermore, \cite[Cor.\ 4.11, 4.12]{AGKS20} gives:
\begin{Theorem}
Let $I\times_f X$ be a generalized cone such that $I$ is open and $X$ is either a geodesic length space that is proper or a locally compact, 
complete length space. Then $I\times_f X$ is globally hyperbolic.
\end{Theorem}
Generalized cones are automatically equipped with the continuous time function $\tau\equiv t: (t,x) \mapsto t$ (see \cite[Lemma 4.2]{AGKS20}). Also, if $X$ is path connected
$Y=I\times_f X$ is scc. For convenience, in what follows when we write that $X$ is a length space we will always tacitly assume that
$X$ possesses only a single \emph{accessibility component} (cf.\ \cite{BBI}), i.e., that any two points in $X$ can be connected by
a path of finite length. In particular, $X$ is always supposed to be path connected.
Let us denote the null distance on $Y$ corresponding to $\tau=t$ by $\hat d_f$. 
By Proposition \ref{3.14}, $\hat d_f$ is continuous.

The following result shows that $\hat d_f$ is a metric (i.e., definite) if $f$ is uniformly bounded below by a positive constant:
\begin{Proposition}\label{prop:f_lower_bound_d_f_def} 
Let $(X,d)$ be a length space and let $Y=I\times_f X$ be the corresponding generalized cone. Suppose  that for some
$f_{\min}\in \R_{>0}$ we have $f_{\min} \le f(t)$ for all $t\in I$.
Then for $p=(t_p,x_p)$, $q=(t_q,x_q)$ 
\begin{align}
\hat d_f(p,q) = |t(p)&-t(q)| = |t_p-t_q|, \qquad\quad\, q \in J^\pm(p) \label{eq:hatdf_J}\\
f_{\min}\cdot d(x_p,x_q) &\le \hat d_f(p,q)  \hspace{6.8em} q\not\in J^\pm(p).\label{eq:hatdf_non_J}
\end{align}
Consequently, $\hat d_f$ is a metric on $Y$ and the time function $\tau$ is locally anti-Lipschitz.
\end{Proposition}
\begin{proof} We adapt arguments from \cite[Lemma 4.9]{AB19}.
If $q \in J^\pm(p)$ and $\beta$ is a causal curve from $p$ to $q$ then by \eqref{eq:causal_min}
\[
\hat d_f(p,q) = \hat L_t(\beta) = |t(q) - t(p)|.
\]
Now let $q\not\in J^\pm(p)$ and take (according to Lemma \ref{lem:existence_of_pw_causal_curves}) a piecewise causal curve $\beta:[a,b]\to Y$
from $p$ to $q$. Let $a=t_0< t_1 <\dots t_k=b$ be a subdivision such that each $\beta_i:=\beta|_{[t_{i-1},t_i]}$ is causal. Changing orientation
if necessary and using \cite[Lemma 3.13]{AGKS20} we may assume that $\beta_i(t) = (t,\alpha_i(t))$ with $\alpha_i$ a
locally Lipschitz curve in $X$. Since $\beta_i$ is causal, \eqref{eq:tl_null_causal} implies that
\[
 v_{\alpha_i}^2(t)\le \frac{1}{f(t)^2} \le \frac{1}{f_{\min}^2}.
\]
Therefore, for $\alpha: J \to X$ the concatenation of the $\alpha_i$ we obtain that 
\[
d(x_p,x_q)\le L_d(\alpha) = \sum_{i=1}^k L_d(\alpha_i) = \sum_{i=1}^k \int_{t_{i-1}}^{t_i} v_{\alpha_i}(t)\,dt \le \frac{1}{f_{\min}}\sum_{i=1}^k |t_i - t_{i-1}|
= \frac{1}{f_{\min}}\hat L_t(\beta).
\]
Taking the infimum over all curves $\beta$ we obtain that $d(x_p,x_q)\le \frac{1}{f_{\min}}\hat d_f(p,q)$.

It is evident from \eqref{eq:hatdf_J}, \eqref{eq:hatdf_non_J} that $\hat d_f$ is definite. The final claim then follows
from Proposition \ref{prop:hat_d_def_tau_antiLip}.
\end{proof}
\begin{Corollary}\label{cor:hat_d_induces_top_warped}
In addition to the assumptions from Proposition \ref{prop:f_lower_bound_d_f_def} let $X$ be locally compact.
Then $\hat d_f$ induces the $D$-topology on $Y=I\times_f X$.
\end{Corollary}
\begin{proof}
This is immediate from the previous result and Proposition \ref{prop:hat_d_induces_top}.
\end{proof}
The proof of Proposition \ref{prop:f_lower_bound_d_f_def} can be utilized to show definiteness of the null distance on
any generalized cone (generalizing \cite[Lemma 3.22]{SV16}):
\begin{Proposition}\label{prop:t_is_definite}
The null distance $\hat d_f$ on any generalized cone $I\times_f X$ (with $X$ a length space) is definite.
\end{Proposition}
\begin{proof}
Let $p=(t_p,x_p)\ne (t_q,x_q)=q \in I\times_f X$. If $t_p\ne t_q$, then $\hat d_f(p,q)\ge |t_p-t_q|>0$ by Proposition \ref{prop:gen_tf} (i).
On the other hand, if $t_p=t_q$ and $x_p\ne x_q$, then there exists some $\delta>0$ such that (at least) a one-sided interval, say $[t_p,t_p+\delta]$
is contained in $I$. Let $\beta$ be a piecewise causal curve connecting $p$ and $q$. If $\beta$ leaves the $\delta$-strip $[t_p,t_p+\delta]$, then
its null-length is at least $\delta$. On the other hand, if $\beta$ remains in $[t_p,t_p+\delta]$, then setting $c:= \min_{t\in [t_p,t_p+\delta]} f(t) >0$, 
the proof of Proposition \ref{prop:f_lower_bound_d_f_def} shows that $\hat L_t(\beta) \ge c d(x_p,x_q) >0$. Taking the infimum over all such 
$\beta$ we obtain that $\hat d_f(p,q)\ge \min (c d(x_p,x_q), \delta)>0$ also in this case.
\end{proof}
Let us denote the space of all piecewise causal paths in $Y$ by $\hat{\mathcal{A}}$. Then we have:
\begin{Theorem}\label{th:length_structure} Let $(X,d)$ be a locally compact length space, $f: I \to \R^+$ continuous, 
and let $Y=I\times_f X$ be the corresponding 
generalized cone. Suppose that $f$ is uniformly bounded below by a positive constant. 
Then $(Y,\hat{\mathcal{A}},\hat L_t)$ defines a length structure on $Y$ in the sense of \cite[Sec.\ 2.1.1]{BBI}.
Moreover, $\hat d_f$ is the intrinsic metric on $Y$ with respect to this length structure, i.e., $(Y,\hat d_f)$ is a length space.
\end{Theorem}
\begin{proof} Since $\hat d_f$ induces the metric topology on $Y$ by Corollary \ref{cor:hat_d_induces_top_warped},
this follows by an obvious adaptation of the proof of \cite[Th.\ 3.5]{AB19}.
\end{proof}

\begin{Lemma}\label{lem:null_paths} Let $Y=I\times_f X$ be a generalized cone, where $(X,d)$ is a  
length space. Then for any $p, q\in Y$ there exists a piecewise null curve connecting $p$ and $q$
(i.e., a piecewise causal curve whose every segment is null).
\end{Lemma}
\begin{proof}
Let $I=\langle a,b\rangle$ and write $p=(t_p,x_p)$, $q=(t_q,x_q)$, where without loss of generality we assume $t_p\le t_q$. 
We first consider the case where in fact $t_p = t_q \in (a,b)$ and pick $\delta>0$ such that $[t_p,t_p+\delta]\sse I$.
If $x_p=x_q$ there is nothing to do. Otherwise, let $\beta: [0,L] \to X$ be a unit speed curve connecting $x_p$ to $x_q$. 
Consider the initial value problem
\begin{equation}\label{eq:ivp}
\left\{
\begin{array}{l}
\dot \alpha_0(s)\, = f(\alpha_0(s))\\
\alpha_0(0) = t_p.
\end{array}
\right.
\end{equation}
The maximal solution to \eqref{eq:ivp} is given by the inverse of the map $r\mapsto \int_{t_p}^r \frac{1}{f(t)}\,dt$
and is defined on $(a_0,b_0)$, where $a_0 = \int_{t_p}^a \frac{1}{f(t)}\,dt$, $b_0= \int_{t_p}^b \frac{1}{f(t)}\,dt$. 
Similarly, for $s_1\in [0,L]$ we consider the final value problem
\begin{equation}\label{eq:fvp}
\left\{
\begin{array}{l}
\dot \alpha_1(s)\, = -f(\alpha_1(s))\\
\alpha_1(s_1) = t_p.
\end{array}
\right.
\end{equation}
The solution $\alpha_1$ of \eqref{eq:fvp} has an analogous explicit description as $\alpha_0$. Combining these 
formulae for $\alpha_0$ and $\alpha_1$ with the fact that $f$ is bounded below and above on $[t_p,t_p+\delta]$
by positive constants we conclude that for any $s_1$ sufficiently small the following conditions are satisfied:
\begin{itemize}
\item $\alpha_0$ and $\alpha_1$ exist on $[0,s_1]$,
\item there exists some $\bar s\in [0,s_1]$ such that $\alpha_0(\bar s) = \alpha_1(\bar s)$, and
\item both $\alpha_0$ and $\alpha_1$ take values in $[t_p,t_p+\delta]$ for all $s\in [0,s_1]$.
\end{itemize}
Consequently, the broken null curve $\gamma: [0,s_1]\to Y$, $\gamma(s) = (\alpha(s),\beta(s))$, where
\[
\alpha(s) = \left\{
\begin{array}{ll}
\alpha_0(s) & \text{for } s\in [0,\bar s]\\
\alpha_1(s) & \text{for } s\in [\bar s,s_1]
\end{array}
\right.
\]
connects the point $(t_p,x_p)$ to the point $(t_p,\beta(s_1))$. Starting from this new point we can iterate the procedure and by 
choosing for $s_1$ a suitably small fraction of $[0,L]$ we obtain the desired broken null curve connecting
$p$ and $q$ in a finite number of steps.

Suppose now that $t_p<t_q$. By introducing, if necessary, an intermediate point $r$ with $t_p < t_r < t_q$ and considering
$p, r$ and $r, q$ separately, we may assume that at most one of $t_p, t_q$ is a boundary point of $I$, say $t_p$ (and the following argument
works just as well if both $t_p$ and $t_q$ lie in the interior of $I$). In this case the solution to \eqref{eq:ivp} attains all $t$-values 
between $t_p$ and $b$, hence in particular the value $t_q$, say at $s=s'$. Then following the null curve $\gamma=(\alpha_0,\beta)$
until $s'$ we obtain a point $(t_q,z)$, which according to the first part of the proof can itself be connected to $q$ 
by a piecewise null curve.
\end{proof}

Based on Lemma \ref{lem:null_paths} we can also prove the following dual to Proposition \ref{prop:f_lower_bound_d_f_def}:
\begin{Proposition}\label{prop:f_upper_bound_d_f_def} 
Let $(X,d)$ be a length space and let $Y=I\times_f X$ be the corresponding generalized cone. Suppose  that for some
$f_{\max}\in \R_+$ we have $f(t)\le f_{\max}$ for all $t\in I$.
Then for $p=(t_p,x_p)$, $q=(t_q,x_q)$ 
\begin{align}
\hat d_f(p,q) &= |t(p)-t(q)| = |t_p-t_q|, \qquad\quad\, q \in J^\pm(p) \label{eq:hatdf_J_new}\\
 \hat d_f(p,q)  &\le f_{\max}\cdot  d(x_p,x_q) \hspace{6.9em} q\not\in J^\pm(p).\label{eq:hatdf_non_J_new}
\end{align}
\end{Proposition}
\begin{proof} Equation \eqref{eq:hatdf_J_new} was already established in Proposition \ref{prop:f_lower_bound_d_f_def}. To show
\eqref{eq:hatdf_non_J_new}, let $\eps>0$ and choose a unit speed path $\beta: [0,L]\to X$ connecting $x_p$ to $x_q$ such 
that $L_d(\beta)=L < d(x_p,x_q)+\eps$. Let us suppose, without loss of generality, that $t_p\le t_q$. As in
the proof of Lemma \ref{lem:null_paths} we then first construct a null curve $\gamma=(\alpha,\beta)$ emanating from $x_p$ and 
reaching an endpoint $r$ with $t_r=t_q$. Since $q\not\in J^+(p)$, the explicit description of $J^+(p)$ in
\cite[Cor.\ 3.24]{AGKS20} shows that  $\alpha$ is defined on $[0,L']$ for some $L'<L$, and again as in the proof of Lemma \ref{lem:null_paths} 
we can then extend $\alpha$ to all of $[0,L]$ such that $\gamma=(\alpha,\beta): [0,L]\to Y$ is a piecewise null curve connecting $p$ and $q$. Denoting the break points of $\gamma$ by $s_0,\dots,s_k$ we then have
\[
\hat L_t(\gamma) = \sum_{i=1}^k |\alpha(s_i)-\alpha(s_{i-1})|. 
\]
Here, since $\gamma$ is null and $v_\beta\equiv 1$, we have $|\alpha(s_i)-\alpha(s_{i-1})| \leq \int_{s_{i-1}}^{s_i} |\dot \alpha(s)|\,ds \le f_{\max}|s_i-s_{i-1}|$. Consequently,
\[
\hat d_f(p,q) \le \hat L_t(\gamma) \le f_{\max}\cdot L < f_{\max}(d(x_p,x_q)+\eps),
\]
giving the claim for $\eps\to 0$.
\end{proof}
\begin{remark}\label{rem:induced_metrics}
Combining Propositions \ref{prop:f_lower_bound_d_f_def} and \ref{prop:f_upper_bound_d_f_def}, it follows that if the warping
function $f$ is bounded from above and below by positive constants then the null distance $\hat d_f$ induces on any fiber $\{t_0\}\times X$ of the 
generalized cone $I\times_f X$ a metric that is (bi-Lipschitz) equivalent to the one induced by the original metric $d$
(by $((t_0,x_1),(t_0,x_2))\mapsto d(x_1,x_2)$). In the special case of
pure Lorentzian products ($f\equiv 1$) both $\hat d_1$ and $d$ induce the same metric on any fiber. This generalizes \cite[Prop.\ 3.3.]{SV16}
and \cite[Lemmas 4.4, 4.9]{AB19}.
\end{remark}
We can now use Propositions \ref{prop:f_lower_bound_d_f_def} and \ref{prop:f_upper_bound_d_f_def} to derive a comparison 
between the null distance $\hat d_f$ on $I\times_f X$ and the Lorentzian product $I\times_1 X$. 
Indeed, we have the following generalization of \cite[Prop.\ 4.10]{AB19}:
\begin{Proposition}\label{prop:null_dist_estimates} Let $I$ be an interval,  
$(X,d)$ a length space, with  corresponding generalized cone $Y=I\times_f X$. 
Suppose that there exist positive constants $f_{\min}, f_{\max}\in \R_+$ such that $0< f_{\min} \le f(t)\le f_{\max}$ for all $t\in I$.
Then for all $p,q\in Y$ we have
\begin{equation}\label{eq:null_dist_bounds}
\min(1,f_{\min}) \hat d_1(p,q) \le \hat d_f(p,q) \le \max(1,f_{\max}) \hat d_1(p,q).
\end{equation}
\end{Proposition}
\begin{proof}
Using the above results, this follows by a straightforward adaptation of the proof of \cite[Prop.\ 4.10]{AB19}.
\end{proof}

\begin{Theorem}\label{th:d_t_definite_and_J} Let $I$ be an interval, $(X,d)$ a length space, and $f:I\to \R^+$ continuous. 
Let $\phi: I\to J\sse \R$ be a strictly monotonically increasing bi-Lipschitz homeomorphism 
and denote by $\tau$ the time function $\tau(t,x)=\phi(t)$ on $I\times_f X$. Then  
for any $p,q \in I\times_f X$ the following are equivalent:
\begin{itemize}
\item[(i)] $p\le q$.
\item[(ii)] $\hat d_{\tau}(p,q) = \tau(q) - \tau(p)$.
\item[(iii)] $\hat d_f(p,q) = t(q) - t(p)$.
\end{itemize}
If $X$ is, in addition, geodesic, then (i)--(iii) is further equivalent to $q\in \overline{I^+(p)}$.
\end{Theorem}
\begin{proof} By Proposition \ref{prop:gen_tf}, (i) implies (ii) and (iii). Since (iii) corresponds to the special
case $\phi= \mathrm{id}_I$, we are left with showing that (ii) implies (i). To see this, we adapt an argument from \cite[Lemma 3.24]{SV16}.
Let $p=(t_p,x_p), q=(t_q,x_q)$ and suppose that $\hat d_\tau(p,q) = \tau(q) - \tau(p)$. If $t_p=t_q$ then $\hat d_\tau(p,q)=0$, so $p=q$ and we are
done (noting that an easy adaptation of Proposition \ref{prop:t_is_definite} yields that $\hat d_\tau$ is definite). Thus suppose that $t_p<t_q$. 
Pick $\eps_0>0$ such that the closed relative $\eps_0$-neighborhood $K_{\eps_0}$ in $J$ of the interval $\phi([t_p,t_q])$ is compact.
For any $0<\eps<\eps_0$ let $\beta$ be piecewise causal from $p$ to $q$ with 
$\hat L_\tau(\beta) \le \tau(q)-\tau(p) +\eps$. From Proposition \ref{lem:null_length_basics} (iii) it follows that $\phi\circ \beta\sse K_\eps\sse K_{\eps_0}$. 
Let $\beta=\beta_1\cdots\beta_k$ be a decomposition of $\beta$ into causal bits. By \cite[Cor.\ 3.13]{AGKS20} we may parametrize each
$\beta_i$ such that $(\pm)\beta_i: [t_i,t_i+\delta_i] \to I\times X$, $\beta_i(t) = (t,\sigma_i(t))$. Here, $t_1=t_p$, $t_k+\delta_k = t_q$.
The concatenation $\sigma$ of the corresponding curves $\sigma_i$ connects $x_p$ to $x_q$, and since $\beta_i$ is causal, for almost
every $t\in [t_i,t_i+\delta_i]$ we have $v_{\sigma_i}(t) \le \frac{1}{f(t)}$. Set 
$c:= \max_{t\in K_{\eps_0}} \frac{(\phi^{-1})'(t)}{f(\phi^{-1}(t))}$.
We have
\begin{align*}
L_d(\sigma) &= \sum_{i=1}^k \int_{t_i}^{t_i+\delta_i} v_{\sigma_i}(t)\,dt = \sum_{i=1}^k \int_{\phi(t_i)}^{\phi(t_i+\delta_i)} v_{\sigma_i}(\phi^{-1}(s)) (\phi^{-1})'(s)\,ds\\
&\le \sum_{i=1}^k \int_{\phi(t_i)}^{\phi(t_i+\delta_i)} \frac{(\phi^{-1})'(s)}{f((\phi^{-1})(s))} \,ds
\end{align*}
Note now that $\sum_{i=1}^k (\phi(t_i+\delta_i) - \phi(t_i)) = \phi(t_q) - \phi(t_q) + (\hat L_\tau(\beta) - (\phi(t_q)-\phi(t_p)))$. Decomposing the above integrals in accordance with this identity and recalling the definition of $\hat L_\tau$ we obtain
\begin{align*}
d(x_p,x_q) \le L_d(\sigma) &\le \int_{\phi(t_p)}^{\phi(t_q)} \frac{(\phi^{-1})'(s)}{f((\phi^{-1})(s))}\,ds + c \cdot (\hat L_\tau(\beta) - (\phi(t_q)-\phi(t_p)))\\
&\le \int_{t_p}^{t_q} \frac{1}{f(t)}\,dt +c \eps.
\end{align*}
Letting $\eps \to 0$ we conclude that $d(x_p,x_q) \le \int_{t_p}^{t_q} \frac{1}{f(t)}\,dt$, which by \cite[Cor.\ 3.24]{AGKS20} means that
$p\le q$, giving (i).

For the final claim we need to show that $\overline{I^+(p)} = J^+(p)$. To this end, 
note that by \cite[Cor.\ 3.24]{AGKS20}, $X$ geodesic implies that $J^+(p)$ is closed, so $\overline{I^+(p)} \sse J^+(p)$.
The converse inclusion follows from the explicit description of $I^+(p)$ and $J^+(p)$ in \cite[Prop.\ 3.22]{AGKS20}
and \cite[Cor.\ 3.24]{AGKS20}.
\end{proof}

\begin{remark} (i) In \cite[Th.\ 3.25]{SV16} the authors consider a Lorentzian warped
product manifold $(I\times S,-dt^2+ f^2(t)h)$ with $I$ an open interval, $f:I\to \R^+$ continuous and $(S,h)$ a 
complete Riemannian manifold. By constructing a suitable conformal metric they then show that for 
any smooth function $\phi:I\to \R^+$ with $\phi'>0$
the time function $\tau(t,x):=\phi(t)$ still encodes the causality of the warped product. 
Theorem \ref{th:d_t_definite_and_J} is a generalization of this result to the metric setting. 

(ii) The assumptions on $\phi$ in Theorem \ref{th:d_t_definite_and_J} can be slightly relaxed (cf.\ \cite[Lemma 3.6]{AGKS20}):
It suffices to assume that $\phi$ is absolutely continuous, monotonically increasing and that $\phi^{-1}$ is absolutely continuous 
with locally bounded derivative.
\end{remark}

\subsection{Convergence of generalized cones and curvature bounds}

In \cite{AB19} convergence of generalized cones based on the null distance has been studied in the spacetime setting. Here we make use of the results just obtained to extend the respective theory beyond the manifold level.
Based on Proposition \ref{prop:null_dist_estimates} we first derive the following essential result on the 
convergence of null-distances for uniformly converging warping functions, generalizing \cite[Prop.\ 5.1]{AB19}:
\begin{Theorem}\label{th:uniform_convergence} 
Let $I$ be an interval, $(X,d)$ a length space, and $f_j: I \to \R^+$ ($j\in \N$) a sequence of continuous functions that converge uniformly to some $f: I \to \R^+$.
Suppose that there exists a uniform lower bound $c$, i.e., $0<c \le f_j(t)$ for all $t\in I$ and all $j\in \N$. Then for the null distances of the corresponding generalized cones $I\times_{f_j} X$, $I\times_f X$ we have:
Let $r>0$ and $p_0, q_0 \in I\times X$. Then 
\begin{equation}\label{eq:unif_conv_d_f_n}
\lim_{j\to \infty} \hat d_{f_j}(p,q) = \hat d_f(p,q)
\end{equation}
uniformly on $B_r^{\hat d_f}(p_0) \times B_r^{\hat d_f}(q_0)$.
\end{Theorem}
\begin{proof} This basically follows by a straightforward adaptation of the proof of \cite[Prop.\ 5.1]{AB19}. Although in
the formulation of that result only pointwise convergence is asserted, it indeed implies the stronger claim made here.
Namely, the arguments laid out there show the following: By uniform convergence, the minimum $f_{\min}$ of $f$ on $I$ is
positive. Given any $\eps \in \big(0,\frac{f_{\min}}{4} \big)$, choose $j_0\in \N$ such that
$\|f-f_j\|_{\infty} < \eps$ for all $j\ge j_0$. Then for all $p,q\in Y$ and all $j\ge j_0$ we have:
\begin{equation}\label{eq:loc_uniform_conv}
\hat d_f(p,q) - \eps \Big(1 + \frac{3}{f_{\min}} \hat d_f(p,q) \Big) \le \hat d_{f_j}(p,q)
\le \hat d_f(p,q) + \eps\Big(1 + \frac{8\eps}{f_{\min}} + \frac{8}{f_{\min}} \hat d_f(p,q)  \Big).
\end{equation}
Using this, it suffices to observe that the factors of $\eps$ in \eqref{eq:loc_uniform_conv} are uniformly bounded on any
$\hat d_f$-ball of finite radius.
\end{proof}
\begin{Corollary}\label{cor:uniform_convergence}
In addition to the assumptions of Theorem \ref{th:uniform_convergence}, suppose that $X$ is locally compact and  that
$\mathrm{diam}(I\times_f X,\hat d_f) < \infty$. 
Then the sequence $(I\times_{f_n} X,\hat d_{f_n})$ of metric spaces converges uniformly to $(I\times_f X,\hat d_f)$ in the sense of
\cite[Def.\ 7.1.5]{BBI}, i.e., $\hat d_{f_n} \rightrightarrows \hat d_f$.
\end{Corollary} 
\begin{proof}
Let $F_n =\mathrm{id}: (I\times X,\hat d_{f_n})\to (I\times X,\hat d_f)$. Then $F_n$ is a homeomorphism since by 
Corollary \ref{cor:hat_d_induces_top_warped}, both $\hat d_{f_n}$ and $\hat d_{f}$ induce the same topology as $D$ on $I\times X$.
Finally, $\hat d_{f_n} \rightrightarrows \hat d_f$ by \eqref{eq:unif_conv_d_f_n} and our assumption on the finite diameter
of $(I\times_f X,\hat d_f)$.
\end{proof}
By Proposition \ref{prop:f_upper_bound_d_f_def}, the diameter assumption of Corollary \ref{cor:uniform_convergence} is satisfied if $X$ is 
a length space, $\mathrm{diam}(X)<\infty$, and either $I$ is compact or $I$ is bounded and the $f_n$ are additionally uniformly bounded from above.
\begin{Corollary}
Let $I$ be an interval and $(X,d)$ a length space. Suppose that $f_n:I\to \R^+$ is a sequence of continuous
functions with a positive uniform lower bound that uniformly converges to $f: I\to \R^+$. 
\begin{itemize}
\item[(i)] Let $X$ be proper 
and fix $p_0=(t_0,x_0)\in I\times X$. Then  
$(I\times_{f_n} X,\hat d_{f_n},p_0) \to (I \times_f X,\hat d_f,p_0)$ in the pointed Gromov-Hausdorff sense.
\item[(ii)] If both $I$ and $X$ are compact, then $(I\times_{f_n} X,\hat d_{f_n}) \to (I \times_f X,\hat d_f)$ in the Gromov-Hausdorff sense.
\end{itemize}
\end{Corollary} 
\begin{proof} (i) Let $r>0$ and $\eps>0$, then by \eqref{eq:unif_conv_d_f_n} we can find $N\in \N$ such that for all $n\ge N$ we have (again 
writing $D^d_r$ for closed $r$-balls with respect to a metric $d$) 
\begin{align}
& D^{\hat d_{f_n}}_{(1-\eps)r}(p_0) \sse D^{\hat d_{f}}_{r}(p_0) \sse D^{\hat d_{f_n}}_{(1+\eps)r}(p_0) \label{eq:ball_incl}\\
& |\hat d_f(p,q) - \hat d_{f_n}(p,q)|<\eps \quad \forall p,q\in D^{\hat d_{f}}_{r}(p_0).\label{eq:df_d_d_f_n}
\end{align}
We now define a map $F: D^{\hat d_{f}}_{r}(p_0) \to D^{\hat d_{f_n}}_{r}(p_0)$ as follows:
If $p\in D^{\hat d_{f_n}}_{(1-\eps)r}(p_0)$ we set $F(p):=p$. Otherwise, 
since $(Y,\hat d_{f_n})$ is a length space (see Theorem \ref{th:length_structure}) and using \eqref{eq:ball_incl}, for any
$p\in D^{\hat d_{f}}_{r}(p_0)$ we can pick some $F(p)\in D^{\hat d_{f_n}}_{r}(p_0)$ such that $\hat d_{f_n}(p,F(p))\le \eps$. 
Then $F(D^{\hat d_{f}}_{r}(p_0)) \supseteq D^{\hat d_{f_n}}_{(1-\eps)r}(p_0)$, hence it is an $\eps$-net in $D^{\hat d_{f_n}}_{r}(p_0)$. 
From \eqref{eq:df_d_d_f_n} and the definition of $F$ we conclude that
\[
|\hat d_f(p,q) - \hat d_{f_n}(F(p),F(q))| \le 3\eps
\]
for all $p, q\in D^{\hat d_{f}}_{r}(p_0)$. Consequently, $F$ is a $3\eps$-isometry and thereby 
\[
d_{GH}((D^{\hat d_{f_n}}_{r}(p_0),p_0),(D^{\hat d_{f}}_{r}(p_0),p_0))\le 6\eps.
\]
(ii) This follows from (i) and general properties of Gromov-Hausdorff convergence (\cite[Prop.\ 2.4]{Jan17}). 
Alternatively, we may use Corollary \ref{cor:uniform_convergence}
and the fact that uniform convergence of compact metric spaces implies Gromov-Hausdorff convergence.
\end{proof}
The following is a compactness result for families of generalized cones:
\begin{Theorem}\label{th:compactness}
Let $I$ be a compact interval and $(X,d)$ a compact length space. Suppose that $\mathcal{F}$ is a family of continuous functions
$I\to \R^+$ that is uniformly bounded above: $\exists C>0$: $f(t)\le C$ for all $t\in I$ and all $f\in \mathcal{F}$. 
Then $\mathcal{Y}:=\{(I\times_f X,\hat d_f) \mid f\in \mathcal{F} \}$ is precompact with respect to the Gromov-Hausdorff topology. Thus
any sequence from $\mathcal{Y}$ possesses a subsequence that converges in the Gromov-Hausdorff sense.
\end{Theorem}
\begin{proof} By compactness of $I$ and $X$, for any $\eps>0$ there exists some $N\in \N$ and $\eps$-nets $t_1,\dots,t_N$ in $I$
and $x_1,\dots,x_N$ in $(X,d)$. Given any $p = (t_p,x_p)\in I\times X$ there then exist $i,j\in \{1,\dots,N\}$ such that $|t_p-t_i|<\eps$
and $d(x_p,x_j)<\eps$. By Proposition \ref{prop:f_upper_bound_d_f_def}, for any $f\in \mathcal{F}$ this implies
\[
\hat d_f(p,(t_i,x_j)) \le \eps \max(1,C).
\]
Consequently, $\{(t_i,x_j)\mid 1\le i,j \le N\}$ is an $\eps\cdot \max(1,C)$-net for $(I\times_f X,\hat d_f)$. This shows that
$\mathcal{Y}$ is uniformly totally bounded (cf.\ \cite[Def.\ 7.4.13]{BBI}). The claim then follows from
\cite[Th.\ 7.4.15]{BBI}.
\end{proof}

For the following result, we recall from \cite[Ex.\ 3.31]{AGKS20} that the Minkowski cone $\mathrm{Cone}(X)$ over 
a geodesic length space from \cite[Sec.\ 2]{AGKS20} can equivalently be represented as the generalized
cone $(0,\infty)\times_{\mathrm{id}} X$. Therefore the GH-convergence result established below in particular 
applies to Minkowski cones.
\begin{Proposition}
Let $(X_n,d_{n},p_n)$ be a sequence of pointed proper length spaces that converge to the pointed proper metric space 
$(X,d,p)$ in the pointed Gromov-Hausdorff sense, and let
$I$ be an interval. If all $I\times_{\mathrm{id}} X_n$ have timelike curvature bounded below by $0$, then 
the same is true of $I\times_{\mathrm{id}} X$.
\end{Proposition} 
\begin{proof} $(X,d)$ is a length space by \cite[Prop.\ 2.7]{Jan17}. Also,
by \cite[Th.\ 2.5]{AGKS20}, $I\times_{\mathrm{id}} X_n$ has timelike curvature bounded below by $0$ if and only if
$(X_n,d_n)$ is an Alexandrov space whose curvature is bounded below by $-1$. This metric curvature bound persists 
through the pointed Gromov-Hausdorff limit of the $(X_n,d_n)$ by \cite[Prop.\ 10.7.1]{BBI} (and the discussion
following \cite[Prop.\ 7.4.12]{BBI}, applied to compact balls containing the quadruples), so appealing again to \cite[Th.\ 2.5]{AGKS20} gives the claim.
\end{proof}
\begin{remark} In particular, the conclusion of the Proposition holds if $I$, $(X,d)$ and all $(X_n,d_n)$ are compact and 
$(X_n,d_n) \to (X,d)$ in the Gromov-Hausdorff sense (cf.\ \cite[Cor.\ 2.5]{Jan17}). 
\end{remark}
We have the following result on null-distance Gromov-Hausdorff convergence of pure Lorentzian products:
\begin{Proposition}\label{prop:GH_implies_GH}
Let $I$ be a compact interval and let $(X_n,d_n)$ be a sequence of compact length spaces. Let 
$(X,d)$ be a compact length space such that $(X_n,d_n) \overset{GH}{\longrightarrow} (X,d)$.
Then $Y_n\equiv (I\times_1 X_n,\hat d_{n,1}) \overset{GH}{\longrightarrow} (I\times_1 X,\hat d_{1})\equiv Y$.
\end{Proposition} 
\begin{proof}
We use the characterization of the Gromov-Hausdorff distance via distortions of correspondences (see \cite[Sec.\ 7.3.3]{BBI}).
For any correspondence $\mathfrak{R} \sse X_n\times X$ between $X_n$ and $X$ we define a correspondence  
$\hat{\mathfrak{R}}\ \sse Y_n\times Y$ between $Y_n$ and $Y$ by
\[
\hat{\mathfrak{R}} := \{((t,x_n),(t,x)) \in Y_n\times Y \mid t\in I,\  (x_n,x)\in \mathfrak{R}  \}.
\]
By Propositions \ref{prop:f_lower_bound_d_f_def} and \ref{prop:f_upper_bound_d_f_def} we have
\begin{equation*}
\hat d_1((s,x),(t,y)) = \left\{ 
\begin{array}{ll}
|s-t| & (s,x)\in J^\pm(t,y)\\
d(x,y) & \text{otherwise.} 
\end{array}
\right.
\end{equation*}
Here, $(s,x)\in J^\pm(t,y)$ means $d(x,y)\le |s-t|$, so $\hat d_1((s,x),(t,y)) = \max(d(x,y),|s-t|)$. Analogously, $\hat d_{n,1}((s,x_n),(t,y_n)) = \max(d_n(x_n,y_n),|s-t|)$.
Thus for the distortions of 
$\mathfrak{R}$ and $\hat{\mathfrak{R}}$ we have
 \begin{align*}
&\mathrm{dis}(\hat{\mathfrak{R}}) = \\
& =\sup \{ |\max(d(x,y),|s-t|) - \max(d_n(x_n,y_n),|s-t|) | \mid  ((s,x_n),(s,x)), ((t,y_n),(t,y)) \in \hat{\mathfrak{R}}\} \\
&\le \sup \{ |d_n(x_n,y_n) - d(x,y)| \mid (x_n,x), (y_n,y) \in \mathfrak{R} \} = \mathrm{dis}({\mathfrak{R}}).
\end{align*}
From this, by \cite[Th.\ 7.3.25]{BBI} we obtain
\[
d_{GH}((Y_n,\hat d_{n,1}),(Y,\hat d_1)) \le \frac{1}{2} \inf_{\mathfrak{R}}(\mathrm{dis}(\hat{\mathfrak{R}})) \le
 \frac{1}{2} \inf_{\mathfrak{R}}(\mathrm{dis}(\mathfrak{R})) = d_{GH}((X_n,d_n),(X,d)),
\]
giving the claim.
\end{proof}
\begin{remark} Since $D^{\hat d_1}_r(t,x) = (I\cap [t-r,t+r])\times D^d_r(x)$ (and analogously for closed $\hat d_n$-balls), 
an easy modification of the previous proof shows the analogous result for pointed Gromov-Hausdorff convergence: 
Let $I$ be any interval and let $(X_n,d_n,x_n)$ be proper length spaces that converge to the proper metric space 
$(X,d,x)$ in the pointed Gromov-Hausdorff sense.
Then if $t_n\to t \in I$, $p_n:=(t_n,x_n)$, $p:=(t,x)$, then also $(I\times_{1}X_n,\hat d_{n,1},p_n) \to (I\times_1 X, \hat d_1,p)$
in the pointed Gromov-Hausdorff sense.
\end{remark}
\begin{Theorem}\label{th:GH_convergence_and_curvature}
Let $I$ be a compact interval and let $(X_n,d_n)$ be a sequence of compact length spaces that converge to 
the compact space $(X,d)$ in the Gromov-Hausdorff sense. If the timelike curvature of each of the 
Lorentzian products $I\times_1 X_n$ is non-negative, then the same is true of $I\times_1 X$.
\end{Theorem}
{\bf Remark.} By Proposition \ref{prop:GH_implies_GH} the assumptions here imply that 
$(I\times_1 X_n,\hat d_{n,1}) \overset{GH}{\longrightarrow} (I\times_1 X,\hat d_{1})$.
\begin{proof}
Any compact length space is geodesic by the Hopf-Rinow Theorem (\cite[Prop.\ I.3.7]{BH99}). Hence we 
may apply \cite[Th.\ 5.7]{AGKS20}, to conclude from our assumption and the fact that $I\times_1 {\mathbb M}^2(0)$ has vanishing
timelike curvature that the metric curvature of $(X_n,d_n)$ is bounded below by $0$.
By \cite[Prop.\ 10.7.1]{BBI}, therefore, also the Gromov-Hausdorff limit $(X,d)$ has non-negative curvature.
Since $(X,d)$ is a length space by \cite[Th.\ 7.5.1]{BBI} (as well as geodesic by what was said above), 
\cite[Th.\ 5.7]{AGKS20} now yields the claim.
\end{proof}

We also have the following convergence result on (timelike) curvature bounds of generalized cones. For its formulation recall from \cite{AB08} that a $C^2$-function $f:I\to (0,\infty)$ is called ($-K'$)-concave (convex), if $f''-K'f\leq 0$ ($\geq 0$).

\begin{Theorem}\label{th:GH_convergence_and_curvature_2}
Let $I$ be compact and assume that $f_n:I\to (0,\infty)$ is a sequence of ($-K_n'$)-concave functions 
that converges to $f:I\to (0,\infty)$ in $C^2$ and such that $K_n'\to K'$. 
Let $(X_n,d_n)$ be a sequence of compact length spaces with lower curvature bounds $K_{X_n}\geq K_n:=\sup_{x\in I}(K_n'f_n^2-(f_n')^2)$ converging to the compact space $(X,d)$ in GH. 
Then the generalized cone $Y=I\times_{f} X$ has timelike curvature bounded below by $K'$.
\end{Theorem}

\begin{proof} Uniform convergence of $f_n$ and $f_n'$ together with $K_n'\to K'$ imply that $K_n \to K:=\sup_{x\in I}(K'f^2-(f')^2)$. Thus by 
\cite[Thm.\ 10.7.1]{BBI}, $X$, being the GH-limit of the $X_n$ has curvature bounded below by $K$. 
Moreover, $f$ is ($-K'$)-concave and so by \cite[Cor.\ 5.4]{AGKS20} the generalized cone $Y=I\times_f X$ has timelike curvature bounded below by $K'$.
\end{proof}

\medskip\noindent
{\bf Acknowledgements.}   
This work was supported by FWF-projects P28770 and P33594 of the Austrian Science Fund FWF.


\begin{thebibliography}{00}

\bibitem{HPS20} Ak\'e Hau, L., Cabrera Pacheco, A.\ J., Solis, D.\ A., On the causal hierarchy of Lorentzian length spaces,
 Classical Quantum Gravity 37 (2020), no. 21, 215013, 22 pp.

\bibitem{AB08} Alexander, S., Bishop, R.L., Lorentz and semi-Riemannian spaces with Alexandrov curvature bounds, {\em Comm.\ Anal.\ Geom.} 16 (2008) no.\ 2, 251--282.

\bibitem{AGKS20} Alexander, S., Graf, M., Kunzinger, M., S\"amann, C., Generalized cones as Lorentzian length spaces: Causality, curvature, and singularity theorems. Submitted, 2020, {\tt arxiv 1909.09575}.

\bibitem{AB19} Allen, B., Burtscher, A., Properties of the null distance and spacetime convergence,
Int.\ Math.\ Res.\ Not.\ (published online 2021), {\tt arXiv 1909.04483}.

\bibitem{AGS:05} Ambrosio, L., Gigli, N., Savaré, G., Gradient flows in metric spaces and in the space of probability measures. Lectures in Mathematics, ETH Z\"urich, 2008.

\bibitem{BH99} Bridson, M.\ R., Haefliger, A.,
\emph{Metric Spaces of Non-positive Curvature}.
Springer-Verlag, Berlin, 1999.

\bibitem{BBI} Burago, D., Burago, Y., Ivanov, S. A Course in Metric Geometry. 
Graduate Studies in Mathematics 33,  American Mathematical Society, Providence, RI, 2001.

\bibitem{CGKE18}  Chru\'sciel, P.T., Grant J.D.E, Kunzinger M, and Minguzzi E, Non-regular spacetime geometry, J. Phys. Conf. Ser. {V}olume 968. 2018.

\bibitem{CG} Chru\'sciel, P.T., Grant, J.D.E.: On Lorentzian causality with
continuous metrics, {\em Classical Quantum Gravity}  29 (2012), no.\ 14, 145001, 32 pp.

\bibitem{GKS19} Grant, J.D.E., Kunzinger, M., S\"amann, C., Inextendibility of spacetimes and 
Lorentzian length spaces, Ann.\ Global Anal.\ Geom.\ 55:133-147 (2019). 

\bibitem{H82} Harris, S.~G., A triangle comparison theorem for Lorentz manifolds, Indiana Univ.\ Math.\ J.\ 31(3), 289--308 (1982).


\bibitem{Jan17} Jansen, D., Notes on pointed Gromov-Hausdorff convergence, {\tt arxiv 1703.09595}.

\bibitem{KS18} Kunzinger, M., S\"amann, C., Lorentzian length spaces,
{\em Ann.\ Global Anal.\ Geom.} 54 (2018), no. 3, 399–447. 

\bibitem{Min19} Minguzzi, E.,
Causality theory for closed cone structures with applications. 
\emph{Rev.\ Math.\ Phys.} 31 (2019), no. 5, 1930001.

\bibitem{Min:19b} Minguzzi, E., 
Lorentzian causality theory.
{\em Living Reviews in Relativity}, 22(1):3, 2019.

\bibitem{ON83}  O'Neill, B., {Semi-Riemannian Geometry. With Applications to
Relativity.} Pure and Applied Mathematics 103. Academic Press, New York, 1983.

\bibitem{SV16} Sormani, C., Vega, C., Null distance on a spacetime, 
{\em Classical Quantum Gravity} 33 (2016), no. 8, 085001, 29 pp. 

\bibitem{SW11} Sormani, C., Wenger, S., The intrinsic flat distance between Riemannian manifolds and other integral current spaces, {\em J.\ Differential Geom.\ } 87 (2011), no. 1, 117--199.

\bibitem{V:21} Vega, C., Spacetime distances: an exploration, arXiv:2103.01191[gr-qc] 





\end{thebibliography}
\end{document}